\documentclass[a4paper,11pt]{amsart}

\usepackage{amssymb}  
\usepackage{mathtools}      
\usepackage{mathabx}        
\usepackage[bb=fourier,cal=euler,scr=rsfs]{mathalfa}	
\usepackage{enumitem}       

\usepackage[ansinew]{inputenc}

 \def\RR{{\mathbb R}}  \def\TT{{\mathbb T}}
   
 \def\ZZ{{\mathbb Z}}

    \def\cW{\mathcal{W}}

\newtheorem*{teo*}{Theorem}

\newtheorem{teo}{Theorem}[section]
\newtheorem{conj}{Conjecture}
\newtheorem{quest}{Question}

\newcommand{\bi}{\begin{itemize}}
\newcommand{\ei}{\end{itemize}}

\theoremstyle{definition}

\theoremstyle{remark}

\author[R. Potrie]{Rafael Potrie}
\address{CMAT, Facultad de Ciencias, Universidad de la Rep\'ublica, Uruguay}
\urladdr{www.cmat.edu.uy/$\sim$rpotrie}
\email{rpotrie@cmat.edu.uy}

\title[Robust dynamics, invariant structures and classification]{Robust dynamics, invariant structures and topological classification}

\thanks{The author was partially supported by CSIC group 618. }

\begin{document}

\maketitle

\begin{abstract}
This text is about geometric structures imposed by robust dynamical behaviour. We explain recent results towards the classification of partially hyperbolic systems in dimension 3 using the theory of foliations and its interaction with topology.  We also present recent examples which introduce a challenge in the classification program and we propose some steps to continue this classification. Finally, we give some suggestions on what to do after classification is achieved. 
\bigskip

\noindent
{\bf Keywords:} Differentiable dynamics, Partial hyperbolicity, Foliations, Topology of 3-manifolds. 

\medskip

\noindent {\bf MSC 2010:} Primary,  37D30  Secondary, 37C20, 57R30,
\end{abstract}

\section{Introduction}\label{SectionIntroduction}
A major goal of dynamics is to be able to predict long term behaviour of a system via the knowledge of the rules that govern the way it is transformed. In this context, whatever can be said \emph{a priori} of a system is relevant.  It is important to search for conditions that can be detected by observing the system evolve in a finite amount of time that will result in consequences on the asymptotic behaviour of it. A beautiful example of this interaction is Shub's entropy conjecture (see e.g. \cite{Shub}) which states that the knowledge of how a manifold wraps around itself (action in homology, detectable in just one iterate) is enough to find lower bounds on the complexity (entropy) of the system.   

The unifying theme of this paper is the dynamical implications of invariant geometric structures and the interaction of the latter with topological and geometric structures on the phase space of the dynamical system. Hyperbolicity has been quite successful in the following sense\footnote{We refer the reader to \cite{Bowen,KH,Shub-libro,Franks-libro} for a more complete mathematical and historical account of hyperbolic theory.}: 

\begin{itemize}
\item it is possible (via cone-fields) to detect if a system is hyperbolic (with given constants) using only finitely many iterates, 
\item hyperbolicity is strongly tied with robust dynamical behaviour, namely, structural stability. (Hyperbolic maps are well known to be structurally (or $\Omega$-)stable and the converse direction is also known to hold in the $C^1$-topology.) This says that one may expect hyperbolicity when the dynamical system, up to change of coordinates, is stable under small perturbations,
\item it is possible to obtain very precise dynamical information of hyperbolic systems from the topological, symbolic and ergodic points of view,
\item for globally hyperbolic systems (Anosov systems) and hyperbolic attractors there is a strong (yet incomplete) classification theory from the topological point of view. 
\end{itemize}

As hyperbolic systems fail to describe all systems one is lead to encounter, it is natural to see what of this theory can extend to some weaker settings and weakenings of the notion of hyperbolicity have appeared in many different ways in the literature since the early 70's. The most ubiquitous generalisations of hyperbolicity are the notions of \emph{non-uniform hyperbolicty} and \emph{partial hyperbolicity}. 

The first notion is a relaxation from the point of view of uniformity which forbids to detect this structure with information of finitely many iterates (though there are many important results that detect this property with positive probability in parametric families of dynamics). Non-uniform hyperbolicity is a property of certain invariant measures and provides very strong implications on the dynamics (see e.g. \cite{Sarig} and references therein). Even in the case of non-uniform hyperbolicity of `large measures' (such as volume), this structure does not impose topological restrictions on the manifolds that admit it \cite{DolgopyatPesin}. Let us remark that in the case of surfaces, there are some soft ways (Ruelle inequality) to detect non-uniform hyperbolicity of certain measures and this allows very strong description of dynamics of smooth diffeomorphisms of surfaces (see e.g. the recent \cite{BCS}).   

Partial hyperbolicity, the main object of this article, has some advantages over non-uniform hyperbolicity, though the study of its dynamics is far from being so developed\footnote{There is an exception though, which is the case of surfaces where there is a quite complete understanding of the dynamics implied by a dominated splitting \cite{PS-dom} (see also \cite{GP}).}. A diffeomorphism $f: M \to M$ is said to admit a \emph{dominated splitting} if it admits a $Df$-invariant continuous splitting $TM= E_1 \oplus \ldots \oplus E_\ell$ with $2\leq \ell \leq \dim M$ into non-trivial subbundles and such that there exists $N>0$ so that for any $x\in M$ and unit vectors $v_i \in E_i(x)$, $v_j \in E_j(x)$ with $i<j$ one has that: 

$$ \| Df^N v_i \| < \| Df^N v_j \|. $$

This condition can be detected via cone-fields (and therefore, given the strength, $N$, by knowing finitely many iterates of $f$, see e.g. \cite[Appendix B]{BDV} or \cite[Chapter 2]{CP}). It is easy to see that in some contexts, this structure already imposes constraints on the topology of manifolds that can admit such diffeomorphisms (e.g. in surfaces, the Euler characteristic must vanish).  

A continuous $Df$-invariant bundle $E$ is said to be \emph{uniformly contracted} (resp. \emph{uniformly expanded}) if there exists $N>0$ so that for every unit vector $v\in E$ one has that:

$$ \|Df^N v\| < 1  \ \text{ (resp. } \ \|Df^{-N} v\| < 1 \ \text{ ).} $$

If a diffeomorphism $f: M \to M$ admits a dominated splitting of the form $TM = E_1 \oplus E_2 \oplus E_3$ (where $E_2$ may be trivial) one says that: 

\begin{itemize}
\item $f$ is \emph{Anosov} if $E_1$ is uniformly contracted, $E_2 = \{0\}$ and $E_3$ is uniformly expanded, 
\item $f$ is \emph{strongly partially hyperbolic} if both $E_1$ is uniformly contracted and $E_3$ is uniformly expanded. 
\item $f$ is \emph{partially hyperbolic} if either $E_1$ is uniformly contracted or $E_3$ is uniformly expanded,
\end{itemize}

We typically put all the uniformly contracting bundles together and denote the resulting bundle as $E^s$ (and symmetrically we denote $E^u$ to the sum of all uniformly expanded subbundles). The rest of the bundles are typically called \emph{center bundles} and it is their existence that makes this structure, on the one hand more flexible and ubiquituous, and on the other, harder to understand. 

\medskip 

As already explained, the study of Anosov systems from a topological point of view is still incomplete, still some remarkable progress has been achieved, notably through the work of Franks, Manning and Newhouse (see \cite[Section 3]{HP-survey} for a fast account with references). Anosov flows have received lot of attention, but even in dimension 3 their classification is still far from complete (see e.g. the introductions of \cite{BarbotFenley,BBY}). It may seem hopeless to try to attack a classification of partially hyperbolic systems, even in dimension 3. 

\smallskip

Pujals has proposed to try to classify strongly partially hyperbolic systems by comparing them with Anosov systems. This would yield relevant information in the quest to understand its dynamics even if the understanding of Anosov systems is incomplete. Later we will try to expand on this. The proposal was undertaken in the pioneering works of Bonatti-Wilkinson \cite{BW} and Brin-Burago-Ivanov \cite{BBI,BI} and has spurred several results in the subject (see \cite{CHHU,HP-survey} for recent surveys). 

The progress in this program has been intense and several unexpected features started to appear \cite{BPP,BGP} leading to the recent realisation of new features of strongly partially hyperbolic diffeomorphisms in dimension 3 in \cite{BGHP}. There is still lot of work ahead, but several positive results give hope that a more precise program can be attacked. We will survey these examples as well as the recent positive results in \cite{BFFP} and we refer the reader to \cite{HP-survey} for a more complete account on the classification of strongly partially hyperbolic systems in 3-manifolds with solvable fundamental group.   

\medskip

{\small {\it Acknowledgments:} I thank the generosity and support of great mathematicians who I admire: C. Bonatti, E. Pujals, A. Wilkinson and specially S. Crovisier and M. Sambarino. I was lucky to work with all my collaborators and students, I cannot explicitly acknowledge all, but I am extremely grateful to each of them, what appears here is full of their contributions. A. Passeggi and A. Sambarino have shared the way since the very beginning. Nati and Amalia make it worth it. }

\section{Robust dynamics}
Partially hyperbolic dynamics arose as a natural generalisation of hyperbolicity \cite{HPS} and as a way to deal with some systems arising naturally in other contexts such as frame flows in negatively curved manifolds \cite{BP}. It also provided examples of structurally stable higher dimensional Lie group actions on manifolds. In analogy to the hyperbolic setting, these systems were shown to be robust and enjoy some stability properties, at least in the case where they admit an invariant foliation tangent to the center direction \cite{Berger, HPS}. Also, as in the case of hyperbolic systems, the strong bundles always integrate into invariant foliations. We refer the reader to \cite{BDV,CP} for recent expositions of these properties. 

Structural stability was early conjectured to imply hyperbolicity (\cite{PalisSmale}) and this was famously proven true by R. Ma\~n\'e (\cite{Manhe1,Manhe2}) in the $C^1$-topology. The ideas introduced by Ma\~n\'e involved the use of \emph{dominated splittings} and have since allowed to characterise other robust dynamical behaviour. 

A system is said to be \emph{transitive} if it has a dense orbit. In \cite{Manhe1} it is shown that a $C^1$-robustly\footnote{i.e. its $C^1$ perturbations are transitive.} transitive surface diffeomorphism must be Anosov (in particular, the surface must be the torus!). These ideas were shown to extend to obtain (weaker) geometric structures in a series of  works \cite{DPU,BDP}. 

\begin{teo}[Bonatti-Diaz-Pujals-Ures]\label{t.BDPU} 
Let $f: M \to M$ be a robustly transitive diffeomorphism. Then, $f$ is \emph{volume hyperbolic}\footnote{This is a weaker notion than partial hyperbolicity. It requires both extremal bundles to verify that the differential uniformly expands/contracts the volume in the bundle, see e.g. \cite{CP} for this and more definitions. When the bundle is one-dimensional this implies that the bundle is uniformly expanded or contracted.}. In particular, for $M$ a 3-dimensional manifold, $f$ is partially hyperbolic. 
\end{teo}

 Parallel to the theory of robust transitivity, another theory deeply related to partial hyperbolicity was developed; it is the theory of \emph{stable ergodicity}, we refer the reader to \cite{Wilkinson} for a survey of this theory. The analog of the above theorem also exists: if a diffeomorphism is stably ergodic in dimension 3, then it is partially hyperbolic \cite{BFP}.  These results are sharp, see the examples in \cite{BD,BV}. 

Theorem \ref{t.BDPU} imposes some restrictions on manifolds admitting robustly transitive (or stably ergodic) diffeomorphisms; e.g. it is known that even-dimensional spheres cannot admit a pair of transverse continuous sub-bundles (see e.g. footnote $1$ in \cite{AB}) and other manifolds may also have such obstructions.  However, these obstructions are far from being as sharp as in the case of surfaces.  For instance, the following is still an open problem:

\begin{quest}\label{q-rt} Is there a 3-manifold which does not admit robustly transitive (or stably ergodic) diffeomorphisms? Does the sphere $S^3$ admit robustly transitive (or stably ergodic) diffeomorphisms? Are there robustly transitive diffeomorphisms of $\TT^3$ homotopic to the identity?
\end{quest}

This seems to be a difficult question. It was shown that the 3-sphere cannot admit strongly partially hyperbolic systems \cite{BI}. For endomorphisms in dimension 2 (which can be seen as a toy model for Question \ref{q-rt}) a complete description of topological obstructions for robust transitivity was recently obtained \cite{LR}. In higher dimensions, even obstructions to the existence of Anosov diffeomorphisms are far from well understood  (see \cite{GL}). 

The ultimate goal would be to provide a topological classification of partially hyperbolic systems which would allow to understand finer dynamical properties that they poses and in that way deduce that if a system has some robust property then its dynamics can be precisely understood. 

\section{Examples}
To the present, we know of the following mechanisms to construct examples (we refer the reader to \cite[Section 3]{CP} for a rather large list of examples) 

\begin{itemize}
\item Algebraic and geometric constructions. Including linear automorphism of tori and nilmanifolds and geodesic and frame flows on negative curvature. 
\item Skew-products. 
\item Examples arising from $h$-transversalities, including Ma\~n\'e type examples (\cite{Manhe0}) and small $C^1$-perturbations of partially hyperbolic systems. 
\item Surgery constructions.
\end{itemize}

None of the four mechanisms is completely understood. Algebraic examples are not even completely understood for the construction of Anosov diffeomorphisms and for geometric examples it is not completely clear which manifolds admit metrics of negative curvature. 

Skew products where the base is more hyperbolic than the fiber can be said to be well understood, but when they work the other way around (the expansion and contraction is seen in the fibers), this is just starting to be studied and several exciting examples are starting to appear (see \cite{HHU-noncoherent,FG,GORH}). 

The concept of $h$-transversality was just recently introduced \cite{BGHP} though it appeared implicitly in several examples. 

Surgery constructions are only partly understood for Anosov flows (see \cite{Fr-Wi, HaTh, Goo, Fried, FenleyAnosov, BBY1}). For partially hyperbolic diffeomorphisms, this kind of construction is in its infancy \cite{Gogolev-surgery, BW}. 

\begin{quest} Are there other ways to construct examples?
\end{quest} 

In the rest of this section we extend a bit on the last two kind of examples.

\subsection{h-Transversalities} 
An $h$-transversality between two partially hyperbolic diffeomorphisms $f, g: M \to M$ is a diffeomorphism $h: M\to M$ which verifies that $Dh(E^{u}_f)$ is transverse to $E^s_g \oplus E^c_g$ and $Dh^{-1}(E^s_g)$ is transverse to $E^c_f \oplus E^u_f$. 

The key property of this condition is that if $f$ is $h$-transverse to itself, then $f^n \circ h$ will be partially hyperbolic for large $n$ (see \cite[Section 2]{BGHP}) so it provides a nice way to construct new examples once one gets enough control on the bundles of a partially hyperbolic diffeomorphism.

 Using some detailed study of these bundles for certain Anosov flows (see \cite{BPP,BGP,BZ}) in \cite{BGHP} we were able to construct a large family of new partially hyperbolic examples. Several questions remain (see in particular \cite[Section 1.4]{BGHP}) but we emphasise on some that involve the $h$-transversality itself: 

\begin{quest} 
Describe the set of $h$-transversalities from an Anosov flow to itself. In particular, is it possible to construct $h$-transversalities from an Anosov flow to itself so that $h$ is isotopic to identity but not through $h$-transversalities? 
\end{quest}
If the last question admits a positive answer one could hope to construct partially hyperbolic diffeomorphisms isotopic to identity behaving very differently from an Anosov flows (see the discussion after Theorem \ref{t.hyperbolic3m} and \cite{BFFP}).

Also, one can wonder if this notion may help creating new examples in higher dimensions, this is completely unexplored territory.

\subsection{Surgeries} 
A conjecture attributed to Ghys (see \cite{Dehornoy}) asserts that every transitive Anosov flow can be obtained (up to topological equivalence) from a given one by performing a finite number of simple operations that consist essentially on making finite lifts or quotients and surgeries (of Fried's type \cite{Fried}). 

We can propose a very vague question in the setting of partially hyperbolic diffeomorphisms: 

\begin{quest}
If one adds operations such as $h$-transversalities and some new type of surgeries to partially hyperbolic diffeomorphisms\footnote{One should also change topological equivalence by some `conjugacy modulo centers' which should be weaker than leaf conjugacy to allow non-dynamical coherence, see also Question \ref{q.conj}.}, can one obtain a classification up to this equivalence? 
\end{quest}

In particular, let us mention that we see all skew-products in dimension 3 as equivalent (the surgery is well explained in \cite[Proposition 4.2]{BW}). Also, the product of an Anosov in $\TT^2$ and the identity in the circle can be easily `surgered' to obtain the time one map of a suspension Anosov flow. On the other hand, we believe that partially hyperbolic diffeomorphisms of \emph{DA}-type should not be in the same equivalence class  as skew products (see \cite{Pot-DCDS} for evidence even if the question is not well posed).

\section{Partial hyperbolicity in other contexts} 
Here we give a small glimpse of other contexts on which one encounters partially hyperbolic dynamics. The choice of topics is certainly biased by the author's interests. 

\subsection{Dynamics far from tangencies}
We refer the reader to \cite{CP} for a more complete account on the relations between partial hyperbolicity and dynamics far from homoclinic tangencies and \cite{Crov-ICM} for a survey on its recent progress. We wish to emphasise the following point though: homoclinic tangencies are a semilocal phenomena, so the partial hyperbolicity obtained by being far from homoclinic tangencies only holds on the chain-recurrent set. In this setting, it makes sense to work and try to analyse the classes independently and so many global arguments (i.e. that depend on the topology of the manifold, or the isotopy class of the diffeomorphism) are lost. 

Probably the main remaining open question in this setting is the following: 

\begin{conj}[Bonatti \cite{Bonatti-Survey}] Generic diffeomorphisms far from homoclinic tangencies have finitely many chain-recurrence classes.
\end{conj}  

See \cite{PujalsSambarino,CSY,CPS,Crov-habilitation} for some progress in this direction. A natural object of study that may combine well global and geometric arguments with the semi-local ones is the study of \emph{attractors}, see \cite{CrPu,CPoS}.

\subsection{Skew-products}

Skew products appear everywhere, as iterated function systems, as so called \emph{fast-slow dynamics} or even as random perturbations of dynamics. The idea is to couple some dynamics with a random, or chaotic behaviour in the base which typically can be modelled by a hyperbolic system. This way, one naturally obtains a partially hyperbolic system (when the coupling is `more random' that the dynamics on the fibres). 

This point of view appears several places in the literature, for example in the notion of \emph{fiber bunching} introduced in \cite{BGMV} by extending ideas of \cite{Led} (see also \cite{AV}). Fiber bunching allows one to see the fibered dynamics as a partially hyperbolic system and construct invariant holonomies which are essentially lifts of strong stable and unstable manifolds to the fibered dynamics. 

We have used this idea to give a somewhat different approach to the study of the Liv\v{s}ic problem for non-conmutative groups \cite{KP} (see \cite{Kalinin} and references therein for an introduction to the problem). This was later used also in recent results such as \cite{Hurtado, AKL}.

\subsection{Discrete subgroups of Lie groups}
We refer the reader to \cite{Benoist-notes,BCLS,GGKW} and reference therein for a more detailed presentation of the subject. We just mention here that in \cite{BPS} we re-interpreted an interesting family of representations of certain groups into Lie groups, known as \emph{Anosov representations} and introduced by Labourie (and later extended to general word-hyperbolic groups by Guichard and Wienhard) in terms of \emph{dominated splitting} and partial hyperbolicity. We refer the reader to \cite{BPS} but we pose here the following question which we believe to be in the same spirit as Theorem \ref{t.BDPU}: 

\begin{quest}
Let $\Gamma$ be a word hyperbolic group, $G$ a semisimple Lie group and $\rho: \Gamma \to G$ a representations which is robustly\footnote{The topology in the space of representations is given by pointwise convergence.} faithful and discrete.  Is it Anosov for some parabolic of $G$? The same question makes sense if one demands $\rho$ to be robustly quasi-isometric and $G$ of real rank $\geq 2$. 
\end{quest}

The question is open even for robustly quasi-isometric representations of the free group in two generators into $\mathrm{SL}(3,\RR)$. This specific question can be posed in the language of linear cocycles as follows:

\begin{quest}
Let $A_0,B_0 \in \mathrm{SL}(3,\RR)$ be two matrices such that for every $A,B$ close to $A_0,B_0$ one has that the linear cocycle over the subshift of $\{A, B, A^{-1}, B^{-1}\}^{\ZZ}$ which does not allow products $AA^{-1}$, $A^{-1}A$, $BB^{-1}$, $B^{-1}B$ verifies that it has positive Lyapunov exponents for every invariant measure. Does this linear cocycle admit a partially hyperbolic splitting? 
\end{quest}

As in the case of Theorem \ref{t.BDPU} this question can be divided in two: show that periodic orbits are \emph{uniformly hyperbolic} at the period, and show that this is enough to obtain a dominated splitting. The second part is \cite[Question 4.10]{BPS} and we were recently able to solve it \cite{KaP}. 

\section{Strong partial hyperbolicity in 3-manifolds} 
In this section we will be concerned with diffeomorphisms $f: M \to M$ where $M$ is a 3-dimensional closed manifold and $f$ is a (strong) partially hyperbolic diffeomorphism admitting an invariant splitting of the form $TM =E^s \oplus E^c \oplus E^u$ into one-dimensional bundles. For some of the progress in higher dimensions we refer the reader to \cite[Section 14]{HP-survey}.

The topological study of these systems can be divided intro three main problems: 
\begin{itemize}
\item to find topological obstructions for $M$ to admit such diffeomorphisms, 
\item to study the integrability of the center bundle, 
\item to classify these systems up to what happens in the center direction.
\end{itemize} 

We present in this section the state of the art in these problems. It is relevant to remark that in dimension 3 there is a quite advanced knowledge on the topology of closed manifolds (see e.g.\cite{Hatcher}) and its interactions with geometry and foliations (we refer the reader to \cite{Calegari} for a nice account).  This allows to pursue a one-by-one methodology to deal with classes of manifolds with increasing complexity. In higher dimensions, different approaches need to be explored and one possibility is to start by making assumptions on the center foliations (in the spirit of \cite{BW}) to work from there instead of studying specific classes of manifolds. 

\subsection{Topological obstructions} 
In this section we are interested with the following:

\begin{quest}
Which 3-manifolds support partially hyperbolic diffeomorphisms? If $M$ admits a partially hyperbolic diffeomorphism, which isotopy classes of diffeomorphisms of $M$ admit partially hyperbolic representatives? 
\end{quest}

This question has a complete answer for 3-manifolds whose fundamental group has subexponential growth. In \cite{BI} it is shown that if the fundamental group of $M$ is abelian, then the action on the homology of $M$ has to be partially hyperbolic (in \cite[Appendix A]{Pot-JMD} it is shown that it has to be strongly partially hyperbolic). This in particular gives that manifolds such as $S^3$ or $S^2 \times S^1$ do not support partially hyperbolic diffeomorphisms. This was extended in \cite{Parwani} to manifolds with subexponential growth of fundamental group obtaining a similar result. By now, we have a complete classification of partially hyperbolic diffeomorphisms in 3-manifolds with (virtually) solvable fundamental group (see \cite{HP-Nil,HP-Sol,HP-survey}) and we know exactly which isotopy classes admit partially hyperbolic representatives. 

On the other hand, in all generality, we do not even know which manifolds admit Anosov flows. For example, the following question is open \footnote{The classification of Anosov flows is completely open in hyperbolic 3-manifolds; when the manifold is toroidal there has been important recent progress towards classification \cite{BarbotFenley1, BarbotFenley2, BarbotFenley, BBY1,BBY}.}:

\begin{quest}
Does every hyperbolic 3-manifold $M$ admit a finite lift supporting an Anosov flow? 
\end{quest}

However, one can pose the following question in the spirit of Pujals' conjecture:

\begin{quest}\label{q.phthenanosov}
If a 3-manifold $M$ with exponential growth of fundamental group admits a (transitive) partially hyperbolic diffeomorphism, does it also (after finite lift) admit a (topological)\footnote{We remark that topological Anosov flows are conjectured to be orbit equivalent to true Anosov flows (\cite{BW}). } Anosov flow? 
\end{quest}

This question is still very far from being understood, but substantial progress has been made in the case where one assumes that $f$ is homotopic to identity. I expect the answer to be yes at least in that case. 

When $M$ is Seifert then the answer to Question \ref{q.phthenanosov} is affirmative (see \cite{HaPS}). On the other hand, even in Seifert manifolds, recent examples \cite{BGP, BGHP} show that several isotopy classes can be produced, but also one can find some obstructions (see \cite[Section 3.1]{BGHP}). 

When $M$ is toroidal and admits Anosov flows transverse to some tori, also several classes of examples can be constructed, and the question of which isotopy classes admit partially hyperbolic representatives is quite open (see \cite{BPP,BGP,BZ,BGHP}).

\subsection{Integrability} 
The stable and unstable bundles of $f$ are known to be (uniquely) integrable into $f$-invariant foliations $\cW^s$ and $\cW^u$. This is for dynamical reasons (see \cite{HPS, HP-survey}). 

However, the lack of regularity of $E^c$ (it is just H\"{o}lder continuous) and the fact that its dynamics is neither contracting or repelling makes its integrability a particular feature. Now, we know several examples where $E^c$ does not integrate into an $f$-invariant foliation. See \cite{HHU-noncoherent, BGHP}. 

Rather than asking for integrability of $E^c$ into an invariant foliation, one typically asks whether $f$ is \emph{dynamically coherent} meaning that both $E^{cs}=E^s\oplus E^c$ and $E^{cu}=E^c \oplus E^u$ integrate into $f$-invariant foliations $\cW^{cs}$ and $\cW^{cu}$. This implies the existence of an $f$-invariant foliation $\cW^c$. This definition involves several subtleties (it is unknown whether the existence of a $f$-invariant center foliation implies dynamical coherence, see \cite{BuW} for discussions on this definition). 

\begin{figure}[ht]
\begin{center}
\includegraphics[scale=0.75]{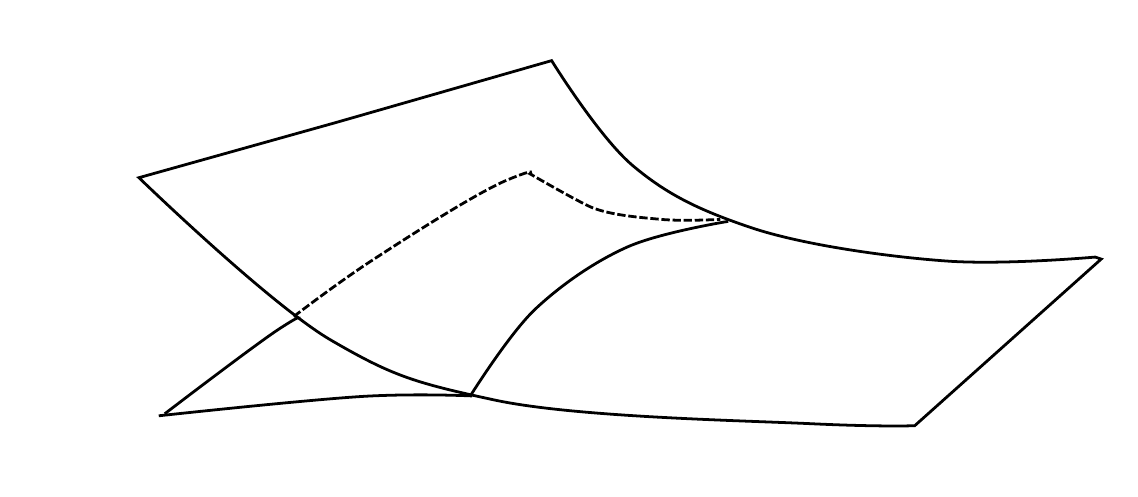}
\begin{picture}(0,0)
\end{picture}
\end{center}
\vspace{-0.5cm}
\caption{{\small Leaves of a branching foliation may merge.}\label{f.branch}}
\end{figure}

A fundamental result was proved by \cite{BI} providing the existence of \emph{branching foliations} (a technical object, see figure \ref{f.branch}) for every strong partially hyperbolic diffeomorphism in 3-manifolds. To many effects, these objects (see also \cite[Chapter 4]{HP-survey}) replace dynamical coherence quite well and allow one to search for classification results. However, it makes sense to ask when do partially hyperbolic diffeomorphisms are dynamically coherent.

We first review some results in this setting: 

\begin{teo}[\cite{BBI2,HHU-noncoherent,Pot-JMD}] 
Let $f: \TT^3 \to \TT^3$ be a partially hyperbolic diffeomorphism. Then, unless there is a torus tangent to $E^{cs}$ or $E^{cu}$ the diffeomorphism $f$ is dynamically coherent.
\end{teo}

This result is proved by obtaining a precise analysis of the structure of the branching foliations provided by \cite{BI} and use this structure to show that branching is not possible (and therefore the branching foliations are indeed foliations).  We would like to emphasise that part of the results obtained in \cite{Pot-JMD} hold for general (not necessarily strong) partially hyperbolic diffeomorphisms through the notion of \emph{almost dynamical coherence}, an open and closed property, which has been also exploited in \cite{FPS, Roldan}.

The same ideas have been pushed into several new contexts by careful combination of topological analysis of the manifolds in hand, comparison to some `model' example and giving some structure to the branching foliations. For instance, with A. Hammerlindl, we were able to show the following results: 

\begin{teo}[\cite{HP-Nil}]
If $N$ is a non-toral nilmanifold and $f: N \to N$ is partially hyperbolic, then $f$ is dynamically coherent. 
\end{teo} 

 \begin{teo}[\cite{HP-Sol}]
 If $M$ is a 3-manifold with (virtually) solvable fundamental group and $f: M \to M$ is partially hyperbolic, then, unless there is a torus tangent $E^{cs}$ or $E^{cu}$ the diffeomorphism $f$ is dynamically coherent. 
 \end{teo}
 
 These results respond affirmatively to a conjecture by Hertz-Hertz-Ures in such manifolds (see \cite{HHU-noncoherent,CHHU}). This conjecture has been disproved recently in \cite{BGHP} but there are still some cases where it can be studied: 
 
 \begin{teo}[\cite{BFFP}] Let $f: M \to M$ be a partially hyperbolic diffeomorphism of a Seifert manifold $M$ such that $f$ is homotopic to identity. Then, $f$ is dynamically coherent. 
 \end{teo}
 
 Seifert manifolds which admit transitive partially hyperbolic diffeomorphisms are, up to finite cover, nilmanifolds or unit tangent bundles of higher genus surfaces (see \cite{HaPS}). In the case of unit tangent bundles of higher genus surfaces, we have constructed in \cite{BGHP} examples which are not dynamically coherent. The configuration responsible for this incoherence is of global nature (see figure \ref{f.nondc}), so it seems natural to ask:

\begin{figure}[ht]
\begin{center}
\includegraphics[scale=0.75]{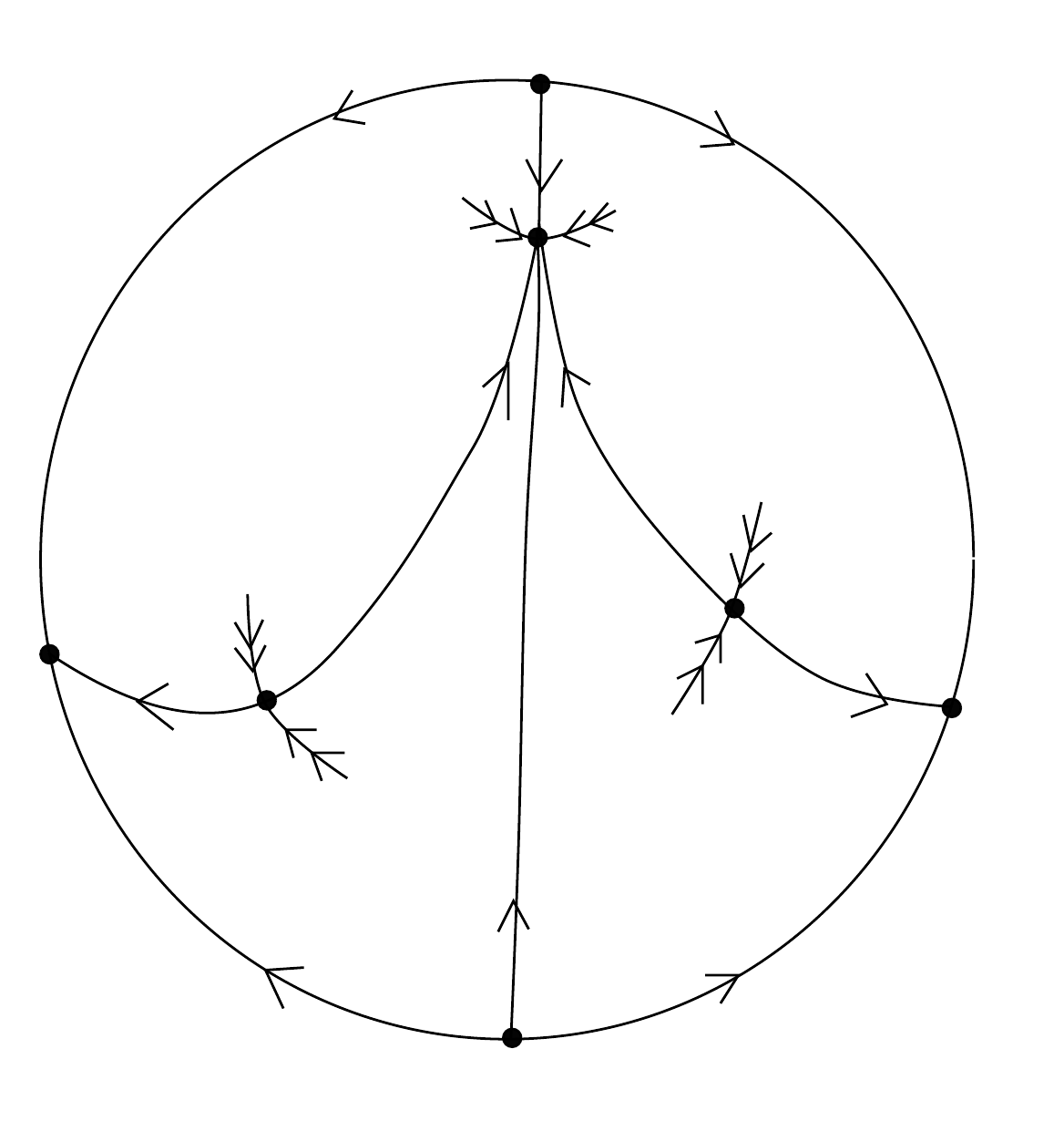}
\begin{picture}(0,0)
\end{picture}
\end{center}
\vspace{-0.5cm}
\caption{{\small The global index in certain periodic $cs$-leaves is negative. Unless $c$-curves merge, there cannot be cancelation so this forces non-dynamical coherence.}\label{f.nondc}}
\end{figure}

 \begin{quest}\label{q-pA}
 For a higher genus surface $S$, if $f: T^1 S \to T^1 S$ is a partially hyperbolic diffeomorphism which induces\footnote{Every diffeomorphism of a Seifert manifold is homotopic to a diffeomorphism preserving the fibers and therefore it has a well defined action on $S$ up to homotopy.} a pseudo-Anosov mapping class in $S$ then $f$ does not admit an $f$-invariant foliation tangent to $E^c$?. 
 \end{quest}
 
 We expect that the techniques used in \cite{BFFP} to classify partially hyperbolic diffeomorphisms of hyperbolic 3-manifolds with $f$-invariant center foliations will allow to provide a positive answer to the previous question. 
 
 When the induced action of $f$ on $S$ is a Dehn-twist (c.f. the examples from \cite{BGP}) then it seems possible that this will imply that the diffeomorphism can be chosen to be dynamically coherent. 
 
 \begin{quest}
 Is there a connected component of partially hyperbolic diffeomorphisms containing both dynamically coherent and non-dynamically coherent diffeomorphisms? 
 \end{quest}
 
 This question is related with the celebrated \emph{plaque expansivity conjecture} (see \cite{HPS, Berger}) since plaque expansivity ensures stability of the invariant foliation, and in case a connected component as in the question exists, it would be natural to check for plaque expansivity in the boundary of dynamically coherent ones.

\subsection{Classification} 
When $f: M \to M$ is a strong partially hyperbolic diffeomorphism which is dynamically coherent, the right notion of classification is given by \emph{leaf conjugacy}: $f,g:M\to M$ dynamically coherent strong partially hyperbolic diffeomorphisms are said to be \emph{leaf conjugate} if there exists a homeomorphism $h: M \to M$ so that $h(\cW^{c}_f(f(x)))= \cW^{c}_g(g(h(x)))$. This notion goes back to \cite{HPS} where a `local stability result' was shown and was retaken in the thesis of A. Hammerlindl \cite{Ham-Thesis} to show that some strong partially hyperbolic diffeomorphisms of $\TT^3$ are leaf conjugate to linear automorphisms of tori. (Notice that this notion does not really require dynamical coherence but the existence of an $f$-invariant center foliation.) 

In this context, we have shown: 

\begin{teo}[\cite{HP-Nil,HP-Sol}]\label{t.HP}
If $M$ is a manifold with (virtually) solvable fundamental group and $f: M \to M$ is a partially hyperbolic diffeomorphism with an $f$-invariant center foliation, then (up to finite lift and iterate) it is leaf conjugate to an algebraic example. 
\end{teo}

In \cite{HP-tori} we further classify those $f$ in such manifolds which do not have an $f$-invariant center foliation. 

In Seifert manifolds, where (topological) Anosov flows are classified (\cite{Ghys,Barbot,Brunella}) one can classify strong partially hyperbolic diffeomorphisms homotopic to identity. 

\begin{teo}[\cite{BFFP}]\label{t.seifert} 
If $f: M \to M$ is a partially hyperbolic diffeomorphism homotopic to identity in a Seifert manifold $M$, then, it is dynamically coherent and leaf conjugate to (up to finite lifts and quotients) the geodesic flow on a surface of negative curvature. 
\end{teo}

We remark that R.Ures has announced a similar result in $T^1S$ which assumes that $f$ is isotopic to the geodesic flow through a path of partially hyperbolic diffeomorphisms (a similar condition to the one studied in \cite{FPS}). 

The results in \cite{BFFP} deal with general partially hyperbolic diffeomorphisms of 3-manifolds which are homotopic to identity. However, there are some points where the precise knowledge of the topology of the manifold under study allows us to give much stronger results. 

For example, when $M$ is a hyperbolic 3-manifold we use the existence of transverse pseudo-Anosov flows to uniform foliations (\cite{Thurston, CalegariPA, Fen2002}) to obtain stronger properties and we deduce the following result (which gives a positive answer to a classification conjecture from \cite{CHHU} for hyperbolic manifolds): 

\begin{teo}[\cite{BFFP}]\label{t.hyperbolic3m}
Let $M$ be a hyperbolic 3-manifold and $f: M \to M$ a dynamically coherent partially hyperbolic diffeomorphism, then $f$ has an iterate which is leaf conjugate\footnote{Technically, the time one map of a topological Anosov flow is not a partially hyperbolic diffeomorphism. But the notion still makes sense.} to the time one map of a (topological) Anosov flow. 
\end{teo}

We make emphasis that this is the first result on classification of partially hyperbolic diffeomorphisms on manifolds where we do not have a model a priori to compare our partially hyperbolic diffeomorphism to. 

This result should hold changing dynamical coherence by the existence of an $f$-invariant foliation tangent to $E^c$. However, I expect that there will be non-dynamically coherent examples in hyperbolic 3-manifolds since some arguments in this theorem resemble those we use in \cite{BGHP} to show that certain examples are not dynamically coherent. These examples would be what we call `double translations' as they act (in the universal cover\footnote{Here, we assume that we take a lift at bounded distance from the identity.}) as translation in both the center-stable and center-unstable (branching) foliations. 
 
Indeed, one can think about the following analogy: Let $f: M \to M$ be a strong partially hyperbolic diffeomorphism  homotopic to the identity on a hyperbolic 3-manifold $M$ which is the suspension of a pseudo-Anosov diffeomorphism $\varphi :S \to S$. Then, the fundamental group of $M$ can be written as a semidirect product $\pi_1(S) \rtimes_{\varphi} \ZZ$. Moreover, the lift of $f$ at bounded distance from the identity in $\tilde M$ commutes with all deck transformations and translates the foliations. Then, we get a group $G_{f}$ which is the direct product of $\langle f\rangle$ and $\pi_1(M)$. 

If one looks at the examples of \cite{BGHP} one has that the fundamental group of $T^1S$ enters in a exact sequence $0 \to \ZZ \to \pi_1(T^1S) \to \pi_1(S) \to 0$ where the inclusion of $\ZZ$ is the center of the group. A lift $\tilde g$ of the diffeomorphism $g: T^1S \to T^1S$ constructed in \cite{BGHP} (which induces a pseudo-Anosov map in the base) will define a group $G_g$ of diffeomorphisms of $\widetilde{T^1S}$ generated by $\pi_1(T^1S)$ and $\tilde g$. Here, the role of $\tilde f$ in $G_f$ is played by the center of $\pi_1(T^1S)$ (which translates center stable and center unstable branching leafs) while $\tilde g$ plays the role of the `semidirect product' in the fundamental group of $M$ above. 

If such an $f$ existed, the dynamics of these two groups in the `circle at infinity' would look quite alike. This suggests on the one hand that examples like this may exist in hyperbolic manifolds, and on the other hand, that it should be possible to answer Question \ref{q-pA} using similar ideas to the ones appearing in Theorem \ref{t.hyperbolic3m}.

\subsection{More questions} 
Several questions remain to be explored in the classification of partially hyperbolic diffeomorphisms in dimension 3. We pose here those we feel are more relevant or that we think might help address the problem of classification. We restrict to the case of 3-manifolds $M$ whose fundamental group is not virtually solvable in view of Theorem \ref{t.HP}. For simplicity we will assume throughout that everything is orientable (the manifold, the bundles, and that $f$ preserves all the orientations). 

Even if a partially hyperbolic diffeomorphism is not dynamically coherent, the one-dimensionality of the center direction allows to integrate the center bundle. However, there may not exist a foliation tangent to it (the strongest integrability that can be ensured beyond the existence of curves tangent to the center is the existence of branching foliations \cite{BI}). 

In the non-dynamically coherent examples of \cite{BGHP} one sees that the space of center leafs is (naturally) homeomorphic to the space of orbits of a geodesic flow in negative curvature (see \cite[Proposition 5.11]{BGHP}) and the dynamics of center leafs is governed by the action of the mapping class of $f$ in the `boundary at infinity'. This could be a general phenomena: 

\begin{quest}\label{q.conj}
Let $f: M \to M$ be a (transitive) partially hyperbolic diffeomorphism. Is there (maybe up to finite cover) an Anosov flow $\phi^t$ on $M$ and continuous degree 1 map $h: M \to M$ sending orbits of $\phi$ to center leafs of $f$?   
\end{quest}

This question is proposing a way to classify dynamics as one could expect that the isotopy class of $f$ will force some dynamics on the center leafs. A related question (maybe more basic, but probably difficult to approach directly)  is: 

\begin{quest}
If $f:M\to M$ is a (transitive) partially hyperbolic diffeomorphism, are all center stable leafs cylinders or planes? 
\end{quest}

Proving this in the case isotopic to identity is an important step towards the proof of Theorems \ref{t.seifert} and \ref{t.hyperbolic3m}, but as far as I am aware this problem has been never attacked directly in general (see \cite{Jinhua} for a positive answer assuming that the dynamics is \emph{neutral} in the center direction). 

The last two questions are definitely related, but they are independent as there are some subtleties in the notions of dynamical coherence, leaf conjugacy, etc. For example I do not know the answer to the following: 

\begin{quest}
If the center foliation of  $f$ is homeomorphic to the orbit foliation of an Anosov flow, are all center stable leaves cylinders or planes?
\end{quest} 

I believe that a reasonable way to attack classification would be to try to understand center leaves at infinity (to be able to avoid taking care of how they merge) and for this the tool of \emph{universal circles} has shown to be quite useful in other contexts (see e.g. \cite{Thurston, Calegari, Fen2002, Frankel}).



\section{Dynamical implications}
We will ignore in this section the very important subject of \emph{conservative} partially hyperbolic diffeomorphisms. They have been extensively treated in other recent surveys such as \cite{CHHU,Wilkinson} with different points of view but great detail. We shall focus mostly on the subject of \emph{robust transitivity} and finitness or uniqueness of attractors for such systems. 

In this direction, one can pose the following question which already appears in \cite{Pot-few} (see also \cite{BGHP}): 

\begin{quest}
Is there an isotopy class of diffeomorphisms of a 3-manifold $M$ such that every strongly partially hyperbolic diffeomorphism in this isotopy class is transitive? Chain-recurrent? 
\end{quest}

In several isotopy classes, such as Anosov times identity on $\TT^3$ the answer is known to be negative (see e.g. \cite{BonattiGuelman,YShi}) for more surprising examples).  But one can still wonder about uniqueness of attractors, or minimal $u$-saturated sets (c.f. \cite{CPoS}). 

Several isotopy classes of (strong) partially hyperbolic diffeomorphisms seem to be now ready for studying subtler dynamical properties. There has been quite some progress in this directions, just to mention a few, we refer the reader to \cite{BV, BFSV, U, CFT, VY} and references therein for advances in the DA-case, \cite{VY1, TY} for the skew-product case and \cite{SaghinYang} for the case of systems leaf conjugate to time one maps of Anosov flows. In all cases it makes sense to try to understand how the entropy behaves\footnote{There is a program in this direction due to J. Buzzi, see \cite{Buzzi}, but it has evolved since.}, how many measures of maximal entropy (or equilibrium states) one may have, physical measures, its statistical properties, its Lyapunov exponents, etc.

I close the paper with a question which I think points towards something we do not really understand yet, and which is more important than the question itself:

\begin{quest}
Let $f: M \to M$ be in the boundary of robustly transitive diffeomorphisms. Does there exist a center Lyapunov exponent for the maximal entropy measure of $f$ which vanishes? 
\end{quest}

It is natural to attack this question in the context of strongly partially hyperbolic diffeomorphisms with one-dimensional center, and it might be that the answer depends on the class. Still, I hope that in the near future we will understand better the transition between the interior of transitive diffeomorphisms and those admitting proper attractors (some progress is in \cite{ACP}, see also \cite[Section 5]{CP}). 





\end{document}